\documentclass[]{elsarticle}

\usepackage{lineno,hyperref}

\usepackage{graphicx}
\usepackage{bm}
\usepackage{multirow}
\usepackage{amsmath}
\usepackage{amsfonts}
\usepackage{amssymb}
\usepackage{amsthm}
\usepackage{secdot}
\usepackage{tabularx}
\usepackage{xcolor,sectsty}

\journal{Computers \& Fluids}

\theoremstyle{plain}

\theoremstyle{definition}

\theoremstyle{remark}

\numberwithin{equation}{section}
\numberwithin{theorem}{section}
\numberwithin{remark}{section}

\newcommand{\Tbold}{\bm{T}}
\newcommand{\nbold}{\bm{n}}

\newcommand{\Ubold}{\bm{U}}

\newcommand{\ubold}{\bm{u}}

\newcommand{\vbold}{\bm{v}}
\newcommand{\wbold}{\bm{w}}

\newcommand{\xbold}{\bm{x}}

\newcommand{\ipt}[2]{\left(#1,#2\right)_{\mathcal{T}_h}}

\newcommand{\ipbt}[2]{\left\langle#1,#2\right\rangle_{\partial \mathcal{T}_h}}

\newcommand{\ipbf}[2]{\left\langle#1,#2\right\rangle_{\mathcal{F}_h}}

\newcommand{\llbracket}{\left[\!\left[}
\newcommand{\rrbracket}{\right] \! \right]}

\newcommand{\llcurve}{\left\{\!\left\{}
\newcommand{\rrcurve}{\right\} \! \right\}}

\DeclareFontFamily{U}{mathx}{\hyphenchar\font45}
\DeclareFontShape{U}{mathx}{m}{n}{
      <5> <6> <7> <8> <9> <10>
      <10.95> <12> <14.4> <17.28> <20.74> <24.88>
      mathx10
      }{}
\DeclareSymbolFont{mathx}{U}{mathx}{m}{n}
\DeclareFontSubstitution{U}{mathx}{m}{n}
\DeclareMathAccent{\widecheck}{0}{mathx}{"71}
\DeclareMathAccent{\wideparen}{0}{mathx}{"75}

\begin{document}

\begin{frontmatter}

\title{Versatile Mixed Methods for Weakly-Compressible Flows}

\author{Edward A. Miller}

\author{David M. Williams \corref{mycorrespondingauthor}} 

\cortext[mycorrespondingauthor]{Corresponding author}
\ead{david.m.williams@psu.edu}

\address{Department of Mechanical Engineering, The Pennsylvania State University, University Park, Pennsylvania 16802}

\begin{abstract}
Versatile mixed finite element methods were originally developed by Chen and Williams for isothermal incompressible flows in ``Versatile mixed methods for the incompressible Navier-Stokes equations," Computers \& Mathematics with Applications, Volume 80, 2020. Thereafter, these methods were extended by Miller, Chen, and Williams to non-isothermal incompressible flows in ``Versatile mixed methods for non-isothermal incompressible flows," Computers \& Mathematics with Applications, Volume 125, 2022. The main advantage of these methods lies in their flexibility. Unlike traditional mixed methods, they retain the divergence terms in the momentum and temperature equations. As a result, the favorable properties of the schemes are maintained even in the presence of non-zero divergence. This makes them an ideal candidate for an extension to compressible flows, in which the divergence does not generally vanish. In the present article, we finally construct the compressible extension of the methods. In addition, we demonstrate the excellent performance of the resulting methods for weakly-compressible flows that arise near the incompressible limit, as well as moderately-compressible flows that arise near Mach 0.5. 
\end{abstract}


\end{frontmatter}


\section{Introduction}
\label{sec;introduction}

In this paper, we discuss the discretization of the compressible Navier-Stokes equations using `versatile mixed finite element methods'. This paper introduces an important extension to previous mixed methods for solving the isothermal incompressible Navier-Stokes equations~\cite{chen2020versatile}, and the \emph{non-isothermal} incompressible Navier-Stokes equations~\cite{miller2022versatile}. 
While mixed finite element methods have seen significant interest for applications to incompressible flows, there has been very limited interest in applying them to compressible flows. This is primarily due to the greater complexity and non-linearity of compressible flows relative to incompressible flows. There are many numerical methods which have the potential to address this challenge. However, in this work we will focus on finite element methods due to their compact stencil, ability to operate on unstructured grids, high-order accuracy, and strong mathematical foundations.

Broadly speaking, we can classify finite element methods for simulating fluids into different categories based on their stabilization strategies. There are at least five stabilization strategies of immediate interest to us: i) residual-based stabilization, ii) numerical-flux-based stabilization, iii) entropy-based stabilization, iv) kinetic-energy-based stabilization, and v) inf-sup stabilization. Often, combinations of these stabilization strategies are used within the same finite element method. In what follows, we will describe each stabilization strategy, its strengths and weaknesses, and then provide some examples of finite element methods which use the strategy. Based on this broad discussion, we will identify the particular stabilization strategies which we deem most effective. Then, we will explain how versatile mixed methods achieve stability using these particular strategies, and how our methods fit into the general landscape of finite element methods for weakly-compressible flows. 

Note: for those readers who are familiar with finite element stabilization strategies, they may skip section~\ref{background_sec}, and proceed directly to section~\ref{current_motivate_sec}.

\subsection{Background} \label{background_sec}

Let us begin by considering a residual-based stabilization strategy. For example, the standard continuous Galerkin (CG) finite element formulation is frequently augmented with a stabilization term which contains the product of the residual operator applied to the solution, a symmetric positive definite matrix, and a residual operator applied to the test functions. Streamline Upwind Petrov-Galerkin (SUPG) methods~\cite{mallet1985finite,hughes1986new,Shakib91} use this strategy for stabilization, as they leverage a convection-based residual operator for the test functions. The SUPG methods are known for their excellent performance in convection-dominated flows, although they are not guaranteed to remain stable in the diffusive limit. Galerkin-Least-Squares (GLS) methods~\cite{Shakib91} are closely related to SUPG methods, with the caveat that they apply the full residual operator (both its advective and diffusive parts) to the test functions. See~\cite{hughes1988computational,hughes1989new,bochev2009least} for a general discussion of least-squares finite element methods. These methods often perform well in both convection- and diffusion-dominated flows, as they are mathematically guaranteed to maintain stability in both settings. Lastly, the Variational Multiscale (VMS) methods~\cite{hughes1998variational,hughes2000large,john2005finite,hughes2007variational} are closely related to GLS methods, with the caveat that they apply the adjoint of the full residual operator to the test functions. In contrast to SUPG and GLS methods, these methods achieve adjoint consistency, although they are not guaranteed to maintain stability. A  clear downside of residual-based stabilization is the substantial difficulty in constructing a well-behaved, positive-definite stabilization matrix. This user-specified matrix must be universal, as it must maintain stability and  recover the correct order of accuracy in smooth parts of the flow, for all applications of interest. 

Next, we consider a numerical-flux-based stabilization strategy. For example, the standard discontinuous Galerkin (DG) finite element formulation is frequently equipped with numerical fluxes which add artificial dissipation to the solution that is proportional to jumps in the solution and/or its gradient. The local DG (LDG)~\cite{nCockburn89a,nCockburn89b,nCockburn90,nCockburn91,nCockburn98}, Bassi-Rebay DG~\cite{Bassi97,bassi1997higha,bassi1997highb}, and compact DG (CDG)~\cite{peraire2008compact,brdar2012compact,pan2022agglomeration} methods use this stabilization strategy. These DG methods have been effectively applied to many convection-dominated flows. Although they are not mathematically guaranteed to maintain stability, the numerical dissipation can be increased as necessary to achieve stabilization. More precisely, the numerical flux contains user-specified constants which can be increased in order to amplify the amount of dissipation. Of course, this strategy is not without risk, as it can result in significant convergence issues within Newton's method due to excessive numerical stiffness. In order to address this issue and to reduce cost, many variants of the classical DG methods have emerged, including hybridized discontinuous Galerkin (HDG)~\cite{cockburn2009unified,fernandez2017hybridized,fernandez2019entropy}, embedded discontinuous Galerkin (EDG)~\cite{cockburn2009analysis,nguyen2015class}, discontinuous Petrov Galerkin (DPG)~\cite{roberts2015discontinuous,ellis2014locally,ellis2016space,rachowicz2021discontinuous}, and variational multiscale discontinuous (VMSD)~\cite{huang2020variational,huang2020adaptive,carson2020provably} methods. 

In addition, we can consider an entropy-based stabilization strategy. For this approach, the compressible Navier-Stokes equations are rewritten in terms of `entropy variables' using a specialized entropy functional~\cite{hughes1986newb}. The development of this functional is based on symmetrizing techniques for conservation laws~\cite{harten1983symmetric,tadmor1984skew}.  Many researchers have applied CG and DG discretization methods to the entropy-symmetrized version of the governing equations, (see the reviews in~\cite{Barth99,Barth06,williams2018entropy,williams2019analysis} for details). The resulting schemes are provably stable for compressible flows. Furthermore, they discretely satisfying the Second Law of Thermodynamics, and can retain stability in the presence of shockwaves at low supersonic speeds without shock-capturing operators. Despite these advantages, there are two significant shortcomings of the entropy-based stabilization approach: a) the entropy variables are highly non-linear functions of the conservative variables, and as a result, they can impede the convergence of Newton's method; b) the entropy is not a useful quantity for controlling the solution as flow conditions approach the incompressible limit. Regarding the last point: the entropy is implicitly related to the compressibility of the flow. For an isothermal, incompressible flow, the entropy remains constant. 


Next, we consider a kinetic-energy-based stabilization strategy. The key idea behind this approach is that kinetic energy is conserved within an unsteady, inviscid, incompressible flow. It turns out that, with the appropriate choice of numerical fluxes and function spaces, this kinetic energy conservation property can be reproduced at the discrete level. An elegant proof of this property for H(div)-conforming mixed methods and specialized DG methods appears in~\cite{guzman2016h}.  It is important to note that these methods are similar (but not identical) in nature to the kinetic energy preserving (KEP) DG methods~\cite{gassner2014kinetic,ortleb2016kinetic,ortleb2017kinetic,flad2017use,klose2020assessing}. These latter methods use a skew-symmetric formulation of the momentum equation in order to create a discrete formulation of the compressible kinetic energy equation which mimics the  continuous equation. The KEP finite element methods are very similar to KEP finite volume methods which were developed earlier in~\cite{jameson2008construction,jameson2008formulation}, and the DG method of~\cite{allaneau2011kinetic}. There are potential advantages for using kinetic-energy-based stabilization instead of entropy-based stabilization, as some KEP schemes possess superior robustness, (see~\cite{jameson2008formulation} for details). Furthermore, kinetic energy is a quantity that remains important in both compressible and incompressible flows, unlike entropy. Lastly, kinetic-energy-based methods avoid the use of highly-non-linear entropy variables. In particular, it is still possible to construct such a method while using conservative or primitive variables. 

Finally, we consider an inf-sup stabilization strategy. For this strategy, the pressure and velocity spaces are required to satisfy an inf-sup compatibility condition.  As a result of this condition, the pressure field remains unique and bounded for incompressible flows, (see~\cite{boffi2013mixed,john2016finite} for details). There are many mixed methods which satisfy the inf-sup condition. For example, H(div)- and H1-conforming mixed methods have been developed for isothermal incompressible flows~\cite{chen2020versatile,schroeder2017pressure,schroeder2018divergence,schroeder2018towards}, and H1-conforming mixed methods have been developed for non-isothermal incompressible flows~\cite{miller2022versatile,Dallmann15}.
A more complete review of inf-sup-stable mixed methods appears in~\cite{chen2020versatile,miller2022versatile}. There has been very limited application of these mixed methods to \emph{compressible} flows. An interesting step in this direction appears in~\cite{schutz2013hybrid,schutz2014adjoint}. Here, a hybrid mixed finite element method is developed for compressible flows which uses a DG method to discretize the convective terms, and an H(div)-conforming method to discretize the diffusive terms. While this approach is quite innovative, it does not utilize the correct pressure and velocity spaces, and thus inf-sup stability is not obtained. In addition, \cite{kweon2003optimal} developed an inf-sup stable mixed method for the stationary, barotropic, linearized, compressible Navier-Stokes equations. Unfortunately, these equations fail to maintain the non-linearity or complexity of the full compressible (or even incompressible) Navier-Stokes equations. A similar argument holds for the mixed methods of~\cite{kellogg1996finite,zhang2011stabilized,egger2018robust}. 


\subsection{Current Work and Motivation} \label{current_motivate_sec}

Our primary goal is to create a mixed finite element method for weakly-compressible flows which uses the most effective and flexible stabilization strategies discussed in the previous section. Towards this end, we have developed versatile mixed methods. These methods utilize three of the five stabilization strategies mentioned above: numerical-flux-based stabilization, kinetic-energy-based stabilization, and inf-sup-based stabilization. We have chosen numerical-flux-based stabilization for our methods instead of residual-based stabilization, due to the inherent simplicity of constructing the numerical fluxes. Next, we have chosen kinetic-energy-based stabilization instead of entropy-based stabilization, due to the overall physical relevance of kinetic energy in both incompressible and compressible environments. Finally, we have chosen inf-sup stabilization due to our interest in developing methods which perform well in nearly-incompressible flows. A direct consequence of the stabilization choices (above), is that versatile mixed methods have very favorable properties. In particular, these methods are provably stable for non-isothermal incompressible flows, as we can rigorously prove L2-stability of the discrete temperature, velocity, and pressure fields. In addition, these methods are accurate, as we can obtain error estimates for all of the discrete fields. The mathematical properties of these methods have been rigorously established in~\cite{chen2020versatile,miller2022versatile}. Furthermore, these properties are maintained in the presence of non-zero divergence of the velocity field. In particular, versatile mixed methods contain divergence terms in the mass, momentum, and temperature equations. The divergence terms in the latter two equations are usually neglected in the construction of conventional mixed methods for incompressible flows. However, these terms are retained in versatile mixed methods. This makes versatile mixed methods particularly suitable for extension to more complex flows which violate the incompressibility constraint. 

It is still necessary for us to actually extend the versatile methods to simulate weakly-compressible flows. This task is the main focus of the present paper.

%
%
With this in mind, it is important for us to briefly discuss other methods which have been developed for weakly-compressible flows. In these flows, numerical methods often encounter problems that stem from discrepancies between wave speeds. This leads to poor conditioning of the matrix system and incorrect predictions of the flow. A way to remedy this issue is the so-called low-Mach-number preconditioners that introduce a preconditioner to the matrix system to bring the wave speeds closer together~\cite{choi1993application,turkel1993review,lee1993progress,weiss1995preconditioning,turkel1999preconditioning}. This approach has seen much success for steady state problems, but by design these methods have a very negative impact on temporal accuracy. A more sophisticated approach has been taken by flux preconditioners, where only the terms in the equations that are naturally dissipative are modified~\cite{guillard1999behaviour,meister2003asymptotic,birken2005stability,birken2005low}. This allows the preservation of the temporal accuracy. In addition, we note that similar preconditioning methods have been developed specifically for finite element methods, (see the work in~\cite{persson2008newton,bassi2009discontinuous,nigro2010discontinuous,klein2016high}). A recent article~(c.f.~\cite{gouasmi2022entropy}) offers an excellent review of research on this topic. While our approach is designed to work for low-Mach-number flows, we do not employ any of the preconditioner approaches outlined above. Instead, our method relies on the usage of function spaces which are normally leveraged to solve the incompressible Navier-Stokes equations. 

We note that similar methods have been developed in recent work, including those described in~\cite{tavelli2017pressure,bermudez2020staggered,zampa2025asymptotic} by Dumbser and colleagues. The schemes in~\cite{tavelli2017pressure,bermudez2020staggered} are pressure-based schemes, and the scheme in~\cite{zampa2025asymptotic} is an entropy-based scheme. In contrast, the present work proposes a temperature-based scheme, which builds on well-known theoretical results for non-isothermal incompressible flows.

\subsection{Overview of the Rest of the Paper}

In section 2, we introduce the governing equations for weakly-compressible flows, and the associated mathematical machinery. In section 3, we outline the versatile mixed methods for weakly-compressible flows. In section 4, we demonstrate that our versatile methods for weakly-compressible flows retain their stability under incompressible conditions, subject to moderate assumptions. In section 5, we present a series of numerical experiments to demonstrate the methods' ability to handle the flows of interest. Lastly, in section 6, we provide a brief summary of our work.


\section{Preliminaries}

Consider the flow of a compressible fluid in a $d$-dimensional domain $\Omega$, where $d = 2$ or 3. Suppose that the fluid has a density field $\rho = \rho \left(t, \xbold \right)$, a momentum field $\rho \ubold = \rho \ubold \left(t, \xbold\right)$, and an internal energy field $\rho e = \rho e \left(t, \xbold\right)$, where $\ubold$ is the velocity and $e$ is the specific internal energy. We assume that $e = C_v T$ where $T$ is the temperature, and $C_v$ is the coefficient of specific heat at constant volume. 

We seek a solution to the motion of the fluid on the time interval $\left(t_0, t_n \right)$ that satisfies the compressible Navier-Stokes equations 
\begin{align}
&\partial_t \, \rho + \nabla \cdot \left( \rho \ubold \right) = S_{\rho}, \label{mass_cons} \\[1.5ex]
&\partial_t \left(\rho \ubold\right) + \nabla \cdot \left( \rho \ubold \otimes \ubold + P \mathbb{I} \right) - \nabla \cdot \left( \rho \bm{\tau} \right)  = \bm{S}_{u}, \label{moment_cons}  \\[1.5ex]
&\partial_t \left(\rho T \right) + \nabla \cdot \left( \rho T  \ubold  \right) - \nabla \cdot \left( \frac{\kappa}{C_v} \nabla T  \right) = - \left(\gamma -1\right) \rho T \left(\nabla \cdot \ubold\right) + \frac{1}{C_v} \left(\rho \bm{\tau} : \nabla \ubold \right)  + S_{T}, \label{energy_cons} 
\end{align}
subject to boundary and initial conditions
\begin{align}
\mathcal{B} \left(\rho, T, \ubold \right) = \bm{0}, \qquad & \text{on} \quad \left[t_0, t_{n} \right] \times \partial \Omega, \label{init_cond} \\[1.5ex]
\rho \left(0,\xbold \right) = \rho_0 \left(\xbold\right), \qquad & \text{in} \qquad \Omega, \\[1.5ex]
T \left(0,\xbold \right) = T_0 \left(\xbold\right), \qquad & \text{in} \qquad \Omega, \\[1.5ex]
\ubold \left(0,\xbold \right) = \ubold_0 \left(\xbold\right), \qquad & \text{in} \qquad \Omega, \label{bound_cond}
\end{align}
where $\partial_t \left( \cdot \right)$ is the temporal derivative operator, $\nabla \left( \cdot \right)$ is the spatial gradient operator, $S_{\rho}$ is a source term for the mass, $\bm{S}_{u}$ is a source term for the linear momentum, $S_T$ is a source term for the internal energy, $P$ is the pressure, $\kappa$ is the coefficient of heat conductivity, $\gamma$ is the ratio of specific heats, and $\mathcal{B} \left(\cdot, \cdot, \cdot \right)$ is a Robin-type boundary condition operator. In addition, $\bm{\tau}$ is the viscous stress tensor
\begin{align}
\bm{\tau} = \nu \left( \nabla \ubold + \nabla \ubold^T - \frac{2}{3} \left(\nabla \cdot \ubold \right) \mathbb{I} \right). \label{stress_tensor}
\end{align}
Here, $\nu = \mu/\rho$ is the kinematic viscosity coefficient, and $\mu$ is the dynamic viscosity coefficient. An explicit formula for $\mu = \mu \left(T\right)$ is given by Sutherland's law
\begin{align}
\mu = \frac{C_{\text{ref}} \, T^{3/2}}{T + S_{\text{ref}}},
\end{align}
where $C_{\text{ref}}$ and $S_{\text{ref}}$ are empirically determined constants. The heat conductivity coefficient $\kappa$ is related to $\mu$ via the following formula
\begin{align}
\kappa = \frac{C_p \, \mu}{\text{Pr}},
\end{align}
where $C_p$ is the coefficient of specific heat at constant pressure and $\text{Pr}$ is the Prandtl number.

In what follows, we will assume that the pressure, temperature, and density are related through an equation of state
\begin{align}
P = \rho R T, \label{ideal_gas}
\end{align}
where $R$ is the specific gas constant. 
%
With this assumption in mind, we can view Eqs.~\eqref{mass_cons} -- \eqref{bound_cond}  as a system of equations for unknowns $\rho$, $T$, and $\ubold$, which is supplemented by the relations in Eqs.~\eqref{stress_tensor} -- \eqref{ideal_gas}. 

Before proceeding further, it is important to note that the system of equations (above) is best-suited for weakly-compressible flows, as Eq.~\eqref{energy_cons} is a transport equation for the temperature, (equivalently, the internal energy). It is often considered best-practice to solve an equation for the \emph{total} energy instead, especially in flows which contain shockwaves. However, in the context of the current paper which focuses on weakly- or moderately-compressible flows, the present  formulation in terms of temperature is sufficient.

We now seek a discrete solution to the continuous equations \eqref{mass_cons}~--~\eqref{bound_cond}. Therefore, we introduce a mesh $\mathcal{T}_h$ of straight-sided, simplex elements, each of which is denoted by $K$. We assume that the domain $\Omega$ is polygonal, and furthermore, that the straight-sided edges of the mesh conform to the geometry of the domain. We choose a mesh in which the individual elements are non-overlapping, and we denote the boundary of each element by $\partial K$. The total collection of faces in the mesh is denoted by $\mathcal{F}_h$, and each face is denoted by $F$. The faces associated with a particular element are denoted by $\mathcal{F}_{K} = \left\{ F \in \mathcal{F}_h : F \subset \partial K \right\}$. The set of all interior faces is denoted by $\mathcal{F}_h^i = \{F \in \mathcal{F}_h : F \cap \partial\Omega = \emptyset \}$, and the set of all boundary faces by $\mathcal{F}_h^{\partial} = \{F \in \mathcal{F}_h : F \cap \partial\Omega \neq \emptyset \}$. We associate each face $F$ with a normal vector $\nbold_F$ which points from the negative (-) side of the face to the positive (+) side. In a similar fashion, the locally defined, outward pointing normal vector for each face of an element is denoted by $\nbold_{F_K}$ or, when the context is clear, simply $\nbold$. 

Next, it is necessary to introduce notation for computing integrals over the elements and faces of the mesh. Suppose that $\phi$ is a scalar function, $\vbold$ and $\wbold$ are vector functions, and $\Tbold$ and $\Ubold$ are tensor functions which are defined on the mesh, and are assumed to be sufficiently smooth. Then, the integrated products of these functions on the mesh are defined as follows
\begin{align*}
\ipt{\phi \, \vbold}{\wbold} & = \sum_{K \in \mathcal{T}_h} \int_{K} \phi \, \vbold \cdot \wbold \, dV, \qquad \ipt{\phi \, \Tbold}{\Ubold} = \sum_{K \in \mathcal{T}_h} \int_{K} \phi \, \Tbold : \Ubold \, dV, \\[1.5ex]
\ipbt{\phi \, \vbold}{\wbold} &= \sum_{K \in \mathcal{T}_h} \int_{\partial K} \phi \, \vbold \cdot \wbold \, dA, \qquad \ipbt{\phi \, \Tbold}{\Ubold} = \sum_{K \in \mathcal{T}_h} \int_{\partial K} \phi \, \Tbold : \Ubold \, dA, \\[1.5ex]
\ipbf{\phi \, \vbold}{\wbold} & = \sum_{F \in \mathcal{F}_h} \int_{F} \phi \, \vbold \cdot \wbold \, dA, \qquad \ipbf{\phi \, \Tbold}{\Ubold} = \sum_{F \in \mathcal{F}_h} \int_{F} \phi \, \Tbold : \Ubold \, dA.
\end{align*}
On a related note, we would like to remind the reader of the following integration by parts formulas 
%
\begin{align*}
\left\langle \phi \vbold, \nbold \right\rangle_{\partial K} &= \left( \phi, \nabla \cdot \vbold \right)_{K} + \left( \vbold, \nabla \phi \right)_{K}, \\[1.5ex]
\left\langle \vbold, \bm{T} \nbold \right\rangle_{\partial K} &= \left(\vbold, \nabla \cdot \bm{T} \right)_{K} + \left(\bm{T}, \nabla \vbold \right)_{K}.
\end{align*}
%
%
Generally speaking, the generic vector function $\vbold$, and the scalar function $\phi$ are not required to be continuous across element boundaries. As a result, it is useful to introduce jump $\llbracket \cdot \rrbracket$ and average $\llcurve \cdot \rrcurve$ operators for the interior faces $F \in \mathcal{F}_h^i$
\begin{align*}
\llbracket \phi \rrbracket &= \phi_{+} - \phi_{-}, \qquad \llbracket \phi \nbold \rrbracket = \phi_{+} \nbold_{+} + \phi_{-} \nbold_{-}, \qquad \llcurve \phi \rrcurve = \frac{1}{2} \left( \phi_{+} + \phi_{-} \right), \\[1.5ex]
\llbracket \vbold \rrbracket &= \vbold_{+} - \vbold_{-}, \qquad \llbracket \vbold \otimes \nbold \rrbracket = \vbold_{+} \otimes \nbold_{+} + \vbold_{-} \otimes \nbold_{-}, \qquad  \llcurve \vbold \rrcurve = \frac{1}{2} \left( \vbold_{+} + \vbold_{-} \right).
\end{align*}
In a similar fashion, for the boundary faces $F \in \mathcal{F}_h^{\partial}$, we define
\begin{align*}
\llbracket \phi \rrbracket &= \phi, \qquad \llbracket \phi \nbold \rrbracket = \phi \nbold, \qquad \llcurve \phi \rrcurve = \phi, \\[1.5ex]
\llbracket \vbold \rrbracket &= \vbold, \qquad \llbracket \vbold \otimes \nbold \rrbracket = \vbold \otimes \nbold, \qquad  \llcurve \vbold \rrcurve = \vbold.
\end{align*}
Next, we introduce convenient function spaces for approximating the density
\begin{align*}
&Q_h^{C} = \left\{ q_h  : q_h \in C^{0} \left( \Omega \right), q_{h} |_{K} \in \mathcal{P}_{k} \left( K \right), \forall K \in \mathcal{T}_h \right\}, \\[1.5ex]
&Q_h^{DC} = \left\{ q_h  : q_h \in L^2 \left( \Omega \right), q_{h} |_{K} \in \mathcal{P}_{k} \left( K \right), \forall K \in \mathcal{T}_h \right\},
\end{align*}
where $\mathcal{P}_{k} \left(K \right)$ is the space of polynomials of degree $\leq k$. One may also approximate the temperature using the function space
\begin{align*}
&R_h^{C} = Q_h^{C} = \left\{ q_h  : q_h \in C^{0} \left( \Omega \right), q_{h} |_{K} \in \mathcal{P}_{k} \left( K \right), \forall K \in \mathcal{T}_h \right\}.
\end{align*}
Finally, one may approximate the velocity field using the Taylor-Hood, Raviart-Thomas, or Brezzi-Douglas-Marini spaces
\begin{align*}
 &\bm{W}_h^{TH} = \left\{ \wbold_h : \wbold_h \in \bm{C}^{0} \left(\Omega\right), \wbold_h |_{K} \in  \left(\mathcal{P}_{k+1} \left( K \right) \right)^d, \forall K \in \mathcal{T}_h \right\}, \\[1.5ex]
&\bm{W}_h^{RT}  = \left\{ \wbold_h : \wbold_h \in \bm{H}\left(\text{div}; \Omega \right), \wbold_h |_{K} \in  \bm{RT}_k \left( K \right), \forall K \in \mathcal{T}_h \right\}, \\[1.5ex]
&\bm{W}_h^{BDM}  = \left\{ \wbold_h : \wbold_h \in \bm{H}\left(\text{div}; \Omega \right), \wbold_h |_{K} \in  \bm{BDM}_{k+1} \left( K \right), \forall K \in \mathcal{T}_h \right\},
\end{align*}
where $\bm{C}^{0} \left(\Omega\right) = \left(C^{0} \left(\Omega\right)\right)^d$ is the vector-valued space of continuous functions, $\bm{RT}_k \left( K \right)$ is the Raviart-Thomas space of degree $k$
\begin{align*}
\bm{RT}_k \left(K \right) = \left(\mathcal{P}_k \left(K \right) \right)^d \oplus \mathcal{P}_k \left(K \right) \xbold,
\end{align*}
and $\bm{BDM}_{k+1} \left(K \right)$ is the Brezzi-Douglas-Marini space of degree $k+1$, whose definition appears in~\cite{boffi2013mixed}.


\section{Extension of Versatile Mixed Methods} \label{mixed_method_sec}

In this section, we define a new class of mixed methods for discretizing Eqs.~\eqref{mass_cons} -- \eqref{energy_cons}. The full derivation of these methods appears in~\ref{method_deriv}. The formal statement of the methods is as follows: 1) consider function spaces $Q_h \subset L^2 \left(\Omega\right)$, $R_h \subset H^1 \left(\Omega\right)$, and $\bm{W}_h \subset \bm{H} \left(\text{div}; \Omega\right)$; 2) choose a set of test functions $\left(q_h, r_h, \wbold_h \right)$ that span $Q_h \times R_h \times \bm{W}_h$; and 3) find $\left(\rho_h, T_h, \ubold_h \right)$ in $Q_h  \times R_h \times \bm{W}_h$ that satisfy: 
\\
\\
\emph{Discrete Mass Equation}
\begin{align}
& \ipt{\partial_t \, \rho_h}{q_h} + \ipt{ \nabla_h \cdot \left(\rho_h  \ubold_h \right)}{q_h} = \ipt{S_{\rho}}{q_h}, \label{mass_cons_disc}
\end{align}
\\
\emph{Discrete Momentum Equation}
\begin{align}
\nonumber & \ipt{\partial_t \left( \rho_h \ubold_h \right)}{\wbold_h} - \ipt{ \left(\rho_h \ubold_h\right) \otimes \ubold_h}{\nabla_h \wbold_h} - R \ipt{\rho_h T_h}{\nabla \cdot \wbold_h} + \ipbt{\widehat{\bm{\sigma}}_{\text{inv}} \, \nbold}{\wbold_h}  \\[1.5ex] 
\nonumber & + \ipt{\rho_h \bm{\tau}_h}{\nabla_h \wbold_h} - \ipbt{\widehat{\bm{\sigma}}_{\text{vis}} \, \nbold}{\wbold_h} \\[1.5ex]
\nonumber & + \ipbt{\widehat{\bm{\varphi}}_{\text{vis}} - \mu_h \ubold_h}{\left( \nabla_h \wbold_h + \nabla_h \wbold_{h}^{T} - \frac{2}{3} \left(\nabla \cdot \wbold_h \right) \mathbb{I} \right) \nbold} \\[1.5ex]
&- \frac{1}{2} \ipt{\Big( \partial_t \, \rho_h + \nabla_h \cdot \left(\rho_h \ubold_h \right) - S_{\rho} \Big) \ubold_h}{\wbold_h}= \ipt{\bm{S}_{u}}{\wbold_h}, \label{moment_cons_disc} 
\end{align}
\\
\emph{Discrete Temperature Equation}
\begin{align}
& \nonumber \ipt{\partial_t \left(\rho_h T_h\right)}{r_h} - \ipt{\rho_h T_h \ubold_h}{\nabla_h r_h} + \ipbt{\widehat{\bm{\phi}}_{\text{inv}}  \cdot \nbold}{r_h}  \\[1.5ex]
\nonumber & + \ipt{\frac{\kappa_h}{C_v} \nabla_h T_h}{\nabla_h r_h}  -\ipbt{\widehat{\bm{\phi}}_{\text{vis}}  \cdot \nbold}{r_h} + \ipbt{\widehat{\lambda}_{\text{vis}} - \frac{\kappa_h}{C_v} T_h}{\nabla_h r_h \cdot \nbold} \\[1.5ex]
\nonumber & -\frac{1}{2} \ipt{\Big(\partial_t \, \rho_h + \nabla_h \cdot \left(\rho_h \ubold_h\right) - S_{\rho} \Big) T_h}{r_h} \\[1.5ex]
\nonumber & = - \left(\gamma-1\right) \Big[ \ipt{\rho_h T_h \left(\nabla \cdot \ubold_h \right)}{r_h} +C_{\text{mod}}\ipt{\left| \rho_h (\nabla\cdot \ubold_h) \right| \nabla_h T_h}{\nabla_h r_h} \Big] \\[1.5ex]
&+ \frac{1}{C_v} \ipt{\rho_h \bm{\tau_h}:\nabla_h \ubold_h}{r_h} + \ipt{S_{T}}{r_h}, \label{energy_cons_disc}
\end{align}
where the quantities with hats (e.g.~$\widehat{\bm{\sigma}}_{\text{inv}}$) denote user-defined numerical fluxes. In addition, we note that $\mu_h = \mu \left(T_h \right)$, $\bm{\tau}_h = \bm{\tau} \left( \rho_h, \mu_h, \ubold_h \right)$, and $\kappa_h = \kappa \left(T_h \right)$. For weakly-compressible flows, we recommend that the numerical fluxes are defined as follows
\begin{align}
\widehat{\bm{\sigma}}_{\text{inv}} &= \llcurve \rho_h \ubold_h \rrcurve \otimes \llcurve \ubold_h \rrcurve + R \llcurve \rho_h T_h \rrcurve \mathbb{I} + \zeta \llcurve \rho_h \rrcurve \left| \ubold_h \cdot \nbold_F \right| \llbracket \ubold_h \otimes \nbold \rrbracket, \label{num_four} \\[1.5ex]
\widehat{\bm{\sigma}}_{\text{vis}} & = \llcurve \rho_h \bm{\tau}_h \rrcurve -\frac{\eta}{h_F} \llcurve \mu_h \rrcurve \llbracket \ubold_h \otimes \nbold \rrbracket, \label{num_five} \\[1.5ex]
\widehat{\bm{\phi}}_{\text{inv}} &= \llcurve \rho_h T_h \rrcurve \ubold_h + \delta \llcurve \rho_h \rrcurve \left| \ubold_h \cdot \nbold_F \right| \llbracket T_h \, \nbold \rrbracket, \label{num_six} \\[1.5ex]
\widehat{\bm{\phi}}_{\text{vis}} & = \frac{1}{C_v} \left( \llcurve  \kappa_h \nabla_h T_h \rrcurve - \frac{\varepsilon}{h_F} \llcurve \kappa_h \rrcurve \llbracket T_h \, \nbold \rrbracket \right), \label{num_seven} \\[1.5ex]
\widehat{\bm{\varphi}}_{\text{vis}} &= \llcurve \mu_h \ubold_h \rrcurve, \qquad \widehat{\lambda}_{\text{vis}} = \frac{1}{C_v} \llcurve \kappa_h T_h \rrcurve,
 \label{num_two}
\end{align}
where $\zeta$, $\eta$, $\delta$, and $\varepsilon$ are adjustable parameters that control the amount of dissipation that is added to the scheme.

We note that the numerical fluxes ensure that conservation is maintained in the discrete momentum and temperature equations. However, the discrete mass equation does not contain a numerical flux. As a result, mass is conserved exactly only if the normal component of the mass flux ($\rho_h \bm{u}_h \cdot \bm{n}$) is continuous at the interfaces between elements. This holds true for a Taylor-Hood-based formulation in which $Q_h \subset H^{1}(\Omega)$, $R_h \subset H^1 \left(\Omega\right)$, and $\bm{W}_h \subset \bm{H} \left(\text{div}; \Omega\right)$. In this case, the product of density and velocity is continuous at the interfaces between elements, even when the density varies. However, if we choose a BDM-based formulation in which $Q_h \subset L^{2}(\Omega)$, $R_h \subset H^1 \left(\Omega\right)$, and $\bm{W}_h \subset \bm{H} \left(\text{div}; \Omega\right)$, then the scheme is only guaranteed to be asymptotically conservative when the density varies.

Returning our attention to the momentum and temperature equations above, one  may observe that we have augmented the schemes by adding `strong residual' terms to the left hand sides of Eqs.~\eqref{moment_cons_disc} and \eqref{energy_cons_disc}
\begin{align*}
   &- \frac{1}{2} \ipt{\Big( \partial_t \, \rho_h + \nabla_h \cdot \left(\rho_h \ubold_h \right) - S_{\rho} \Big) \ubold_h}{\wbold_h}, \\[1.5ex]
   &-\frac{1}{2} \ipt{\Big(\partial_t \, \rho_h + \nabla_h \cdot \left(\rho_h \ubold_h\right) - S_{\rho} \Big) T_h}{r_h}.
\end{align*}
These are skew-symmetrizing terms which maintain consistency, while helping to stabilize the schemes. In particular, they ensure that the convective operators in the momentum and temperature equations become semi-coercive in the incompressible limit. One may consult Lemma 6.4 of~\cite{chen2020versatile} and Lemma~4.12 of~\cite{miller2022versatile} for details. 

Lastly, we have added the following term to the right hand side of Eq.~\eqref{energy_cons_disc}
\begin{align*}
    -C_{\text{mod}} \left(\gamma - 1\right) \ipt{\left| \rho_h (\nabla\cdot \ubold_h) \right| \nabla_h T_h}{\nabla_h r_h},
\end{align*}
where $C_{\mathrm{mod}}$ is a stabilization constant.
This term allows us to control the temperature field in flows which are dominated by temperature-dependent buoyancy effects. It is set to zero in most cases.


\section{Incompressible Stability}

In this section, we introduce a new type of non-linear stability for finite element methods. A finite element method for solving the compressible Navier-Stokes equations is said to possess \emph{incompressible stability} or equivalently is said to be \emph{incompressibly stable} if, upon setting $\rho_h = \text{const}$, $\mu_h = \text{const}$, $\kappa_h = \text{const}$, ${S_{\rho} = 0}$, and neglecting or controlling the viscous dissipation term, we recover an inf-sup stable method for the non-isothermal, \emph{incompressible} Navier-Stokes equations. Broadly speaking, enforcing incompressible stability is a way of enforcing compatibility between the finite element discretizations for compressible and incompressible flows. We contend that, a method for compressible flows which possesses incompressible stability is more likely to maintain robust behavior in the incompressible limit.

 In what follows, we will establish that the finite element methods in Eqs.~\eqref{mass_cons_disc}--\eqref{energy_cons_disc} are incompressibly stable. Towards this end, we initially set $\rho_h = \rho_0 = \text{const}$ in Eq.~\eqref{ideal_gas}
 \begin{align*}
    P_h &= \rho_0 R T_h,
\end{align*}
or equivalently
\begin{align}
    p_h & \equiv R T_h, \label{kinematic_press}
 \end{align}
 where we have defined $p_h = P_h/\rho_0$ as the kinematic pressure. We can immediately observe that $p_h$ and $T_h$ now reside in the same function space, i.e. $(p_h, T_h) \in R_h \times R_h$. Furthermore, if we choose $Q_h = R_h = Q_h^{C}$, then $p_h \in Q_h$. 

 
 We can use the observations above to rewrite Eqs.~\eqref{mass_cons_disc}--\eqref{energy_cons_disc}. Upon performing this operation, and setting $\rho_h = \rho_0 = \text{const}$, $\mu_h = \mu_0 = \text{const}$, $\kappa_h = \kappa_0 = \text{const}$, and $S_{\rho} =0$ in Eqs.~\eqref{mass_cons_disc}--\eqref{energy_cons_disc}, we obtain the following simplified equations:
 \\
 \\
\emph{Discrete Mass Equation}
 \begin{align}
 \ipt{ \nabla_h \cdot \ubold_h}{q_h} = 0, \label{mass_cons_disc_incomp}
\end{align}
\\
\emph{Discrete Momentum Equation}
\begin{align}
\nonumber & \ipt{\partial_t \ubold_h}{\wbold_h} - \ipt{\ubold_h \otimes \ubold_h}{\nabla_h \wbold_h} -  \ipt{p_h}{\nabla \cdot \wbold_h} + \ipbt{\frac{\widehat{\bm{\sigma}}_{\text{inv}}}{\rho_0} \, \nbold}{\wbold_h}  \\[1.5ex] 
\nonumber & + \ipt{\bm{\tau}_h}{\nabla_h \wbold_h} - \ipbt{\frac{\widehat{\bm{\sigma}}_{\text{vis}}}{\rho_0} \, \nbold}{\wbold_h} \\[1.5ex]
\nonumber & + \ipbt{\frac{\widehat{\bm{\varphi}}_{\text{vis}}}{\rho_0} - \nu \ubold_h}{\left( \nabla_h \wbold_h + \nabla_h \wbold_{h}^{T} - \frac{2}{3} \left(\nabla \cdot \wbold_h \right) \mathbb{I} \right) \nbold} \\[1.5ex]
&- \frac{1}{2} \ipt{ \left( \nabla_h \cdot \ubold_h \right)  \ubold_h}{\wbold_h}= \ipt{\bm{f}_{u}}{\wbold_h}, \label{moment_cons_disc_incomp} 
\end{align}
\\
\emph{Discrete Temperature Equation}
\begin{align}
& \nonumber \ipt{\partial_t T_h}{r_h} - \ipt{T_h \ubold_h}{\nabla_h r_h} + \ipbt{\frac{\widehat{\bm{\phi}}_{\text{inv}}}{\rho_0}  \cdot \nbold}{r_h}  \\[1.5ex]
\nonumber & + \gamma \alpha \ipt{\nabla_h T_h}{\nabla_h r_h}  -\ipbt{\frac{\widehat{\bm{\phi}}_{\text{vis}}}{\rho_0}  \cdot \nbold}{r_h} + \ipbt{\frac{\widehat{\lambda}_{\text{vis}}}{\rho_0} - \gamma \alpha T_h}{\nabla_h r_h \cdot \nbold} \\[1.5ex]
\nonumber & -\frac{1}{2} \ipt{ \left( \nabla_h \cdot \ubold_h \right)  T_h}{r_h} \\[1.5ex]
 \nonumber & = - \left(\gamma-1\right) \Big[ \ipt{ \left(\nabla \cdot \ubold_h \right) T_h}{r_h} +C_{\text{mod}}\ipt{\left| \nabla\cdot \ubold_h \right| \nabla_h T_h}{\nabla_h r_h} \Big]  \\[1.5ex]
 &+ \frac{1}{C_v} \ipt{\bm{\tau_h}:\nabla_h \ubold_h}{r_h} + \ipt{f_{T}}{r_h}, \label{energy_cons_disc_incomp}
\end{align}
where we have set $\alpha \equiv \kappa_{0}/C_{p} \rho_0$, $\bm{f}_{u} \equiv \bm{S}_{u}/\rho_0$, and $f_{T} \equiv S_{T}/\rho_0$. Next, upon substituting $\rho_h = \rho_0 = \text{const}$, $\mu_h = \mu_0 = \text{const}$, and $\kappa_h = \kappa_0 =\text{const}$ into the numerical fluxes (Eqs.~\eqref{num_four}--\eqref{num_two}) and dividing through the result by $\rho_0$, we have that
\begin{align}
\widecheck{\bm{\sigma}}_{\text{inv}} \equiv \frac{\widehat{\bm{\sigma}}_{\text{inv}}}{\rho_0} &= \llcurve \ubold_h \rrcurve \otimes \llcurve \ubold_h \rrcurve +  \llcurve p_h \rrcurve \mathbb{I} + \zeta \left| \ubold_h \cdot \nbold_F \right| \llbracket \ubold_h \otimes \nbold \rrbracket, \label{num_four_incomp} \\[1.5ex]
\nu \widecheck{\bm{\sigma}}_{\text{vis}} \equiv \frac{\widehat{\bm{\sigma}}_{\text{vis}}}{\rho_0} & = \llcurve \bm{\tau}_h \rrcurve -\frac{\eta \nu}{h_F} \llbracket \ubold_h \otimes \nbold \rrbracket, \label{num_five_incomp} \\[1.5ex]
\widecheck{\bm{\phi}}_{\text{inv}} \equiv \frac{\widehat{\bm{\phi}}_{\text{inv}}}{\rho_0} &= \llcurve T_h  \rrcurve \ubold_h + \delta \left| \ubold_h \cdot \nbold_F \right| \llbracket T_h \, \nbold \rrbracket, \label{num_six_incomp} \\[1.5ex]
\gamma \alpha \widecheck{\bm{\phi}}_{\text{vis}} \equiv \frac{\widehat{\bm{\phi}}_{\text{vis}}}{\rho_0} & = \gamma \alpha \left( \llcurve  \nabla_h T_h \rrcurve - \frac{\varepsilon}{h_F}  \llbracket T_h \, \nbold \rrbracket \right), \label{num_seven_incomp} \\[1.5ex]
\nu \widecheck{\bm{\varphi}}_{\text{vis}} \equiv \frac{\widehat{\bm{\varphi}}_{\text{vis}}}{\rho_0} &= \nu \llcurve \ubold_h \rrcurve, \qquad \gamma \alpha \widecheck{\lambda}_{\text{vis}} \equiv \frac{\widehat{\lambda}_{\text{vis}}}{\rho_0} = \gamma \alpha \llcurve T_h \rrcurve.
 \label{num_two_incomp}
\end{align}
We may then substitute Eqs.~\eqref{num_four_incomp}--\eqref{num_two_incomp} into Eqs.~\eqref{mass_cons_disc_incomp}--\eqref{energy_cons_disc_incomp} in order to obtain
 \begin{align}
 \ipt{ \nabla \cdot \ubold_h}{q_h} = 0, \label{mass_cons_disc_incomp_simp}
\end{align}
\begin{align}
\nonumber & \ipt{\partial_t \ubold_h}{\wbold_h} - \ipt{\ubold_h \otimes \ubold_h}{\nabla_h \wbold_h} -  \ipt{p_h}{\nabla \cdot \wbold_h} + \ipbt{\widecheck{\bm{\sigma}}_{\text{inv}} \, \nbold}{\wbold_h}  \\[1.5ex] 
\nonumber & + \ipt{\bm{\tau}_h}{\nabla_h \wbold_h} - \nu \ipbt{\widecheck{\bm{\sigma}}_{\text{vis}} \, \nbold}{\wbold_h} \\[1.5ex]
\nonumber & + \nu \ipbt{\widecheck{\bm{\varphi}}_{\text{vis}} - \ubold_h}{\left( \nabla_h \wbold_h + \nabla_h \wbold_{h}^{T} - \frac{2}{3} \left(\nabla \cdot \wbold_h \right) \mathbb{I} \right) \nbold} \\[1.5ex]
&- \frac{1}{2} \ipt{ \left( \nabla \cdot \ubold_h \right)  \ubold_h}{\wbold_h}= \ipt{\bm{f}_{u}}{\wbold_h}, \label{moment_cons_disc_incomp_simp} 
\end{align}
\begin{align}
& \nonumber \ipt{\partial_t T_h}{r_h} - \ipt{T_h \ubold_h}{\nabla_h r_h} + \ipbt{\widecheck{\bm{\phi}}_{\text{inv}}  \cdot \nbold}{r_h}  \\[1.5ex]
\nonumber & + \gamma \alpha \ipt{\nabla_h T_h}{\nabla_h r_h}  - \gamma \alpha \ipbt{\widecheck{\bm{\phi}}_{\text{vis}}  \cdot \nbold}{r_h} + \gamma \alpha \ipbt{\widecheck{\lambda}_{\text{vis}} - T_h}{\nabla_h r_h \cdot \nbold} \\[1.5ex]
\nonumber & -\frac{1}{2} \ipt{ \left( \nabla \cdot \ubold_h \right)  T_h}{r_h} \\[1.5ex]
\nonumber & = - \left(\gamma-1\right) \Big[ \ipt{ \left(\nabla \cdot \ubold_h \right) T_h}{r_h} +C_{\text{mod}}\ipt{\left| \nabla\cdot \ubold_h \right| \nabla_h T_h}{\nabla_h r_h} \Big] \\[1.5ex] 
 &+ \frac{\nu}{C_v} \ipt{\left(\nabla_h \bm{u}_h + \nabla_h \bm{u}_{h}^{T} - \frac{2}{3}\left(\nabla \cdot \bm{u}_h\right) \mathbb{I} \right):\nabla_h \ubold_h}{r_h} + \ipt{f_{T}}{r_h}. \label{energy_cons_disc_incomp_simp}
\end{align}
The resulting class of methods is identical to the class which was originally introduced in section 3 of~\cite{miller2022versatile}, with two important caveats: 
\begin{itemize}
    \item There is a viscous dissipation term
    \begin{align}
        \frac{\nu}{C_v} \ipt{\left(\nabla_h \bm{u}_h + \nabla_h \bm{u}_{h}^{T} - \frac{2}{3}\left(\nabla \cdot \bm{u}_h\right) \mathbb{I} \right):\nabla_h \ubold_h}{r_h}, \label{visc_disp}
    \end{align}
    which now appears on the right hand side of the discrete temperature equation (Eq.~\eqref{energy_cons_disc_incomp_simp}). This term was not present in section 3 of the original work.
    \item The temperature space $R_h$ is one degree lower than in the original work.
\end{itemize}
We note that the original methods in~\cite{miller2022versatile} were proven to be inf-sup stable for the incompressible, non-isothermal, Navier-Stokes equations; and furthermore this fact is unaffected by the change in the degree of the temperature space. However, inf-sup stability was \emph{not} proven in the presence of the viscous dissipation term (Eq.~\eqref{visc_disp}). In fact, our analysis of the discrete temperature equation suggests that this term has a de-stablizing effect. More prescisely, upon substituting $r_h = T_h$ into Eq.~\eqref{energy_cons_disc_incomp_simp}, one obtains
\begin{align}
    \frac{1}{2} \frac{d}{dt} \left\| T_h \right\|_{L^2(\Omega)}^{2} = \text{other terms} + \frac{\nu}{C_v} \ipt{\left(\nabla_h \bm{u}_h + \nabla_h \bm{u}_{h}^{T} - \frac{2}{3}\left(\nabla \cdot \bm{u}_h\right) \mathbb{I} \right):\nabla_h \ubold_h}{T_h}
 \end{align}
 where evidently
 \begin{align}
     \frac{\nu}{C_v} \ipt{\left(\nabla_h \bm{u}_h + \nabla_h \bm{u}_{h}^{T} - \frac{2}{3}\left(\nabla \cdot \bm{u}_h\right) \mathbb{I} \right):\nabla_h \ubold_h}{T_h} \geq 0.
 \end{align}
 Therefore, the viscous dissipation term is a source term, which can cause the $L^2$-norm of $T_h$ to grow. The presence of the viscous dissipation term can be reconciled in one of the following ways:
\begin{enumerate}
    \item One may assume that the viscous dissipation term vanishes in the incompressible limit. This assumption holds true at the continuous level---see the analysis in Appendix A.1 of~\cite{miller2022versatile}. However, this assumption is not guaranteed to hold true at the discrete level. Nevertheless, in section 6 of~\cite{miller2022versatile}, we performed numerical experiments which indicated that the viscous dissipation term has no significant impact in the incompressible flows we tested.
    \item One may modify the compressible temperature equation (Eq.~\eqref{energy_cons_disc}) by adding a Smagorinsky-type dissipation term~\cite{smagorinsky1963general} to the right hand side
    \begin{align}
        -\frac{1}{C_v} \ipt{\rho_h \bm{\tau}_{h}^{\mathrm{turb}}:\nabla_h \ubold_h}{r_h}, 
    \end{align}
    where
    \begin{align}
        \bm{\tau}_{h}^{\mathrm{turb}} &= \nu_{h}^{\mathrm{turb}} \left( \nabla_h \ubold_h + \nabla_h \ubold_h^T - \frac{2}{3} \left(\nabla \cdot \ubold_h \right) \mathbb{I} \right), \\[1.5ex]
        \nu_{h}^{\mathrm{turb}} &= C_{\mathrm{smag}} h^2 \sqrt{ \frac{1}{2} \left( \nabla_h \ubold_h + \nabla_h \ubold_h^T \right) : \left( \nabla_h \ubold_h + \nabla_h \ubold_h^T \right)},
    \end{align}
    and $C_{\mathrm{smag}}$ is a modeling constant.
    Evidently, for sufficiently large values of $C_{\mathrm{smag}}$, we obtain a stable discretization in the incompressible limit, as
\begin{align}
    \frac{1}{2} \frac{d}{dt} \left\| T_h \right\|_{L^2(\Omega)}^{2} = \text{other terms} +  \frac{1}{C_v} \ipt{\left(\nu - \nu_h^{\mathrm{turb}}\right) \left(\nabla_h \bm{u}_h + \nabla_h \bm{u}_{h}^{T} - \frac{2}{3}\left(\nabla \cdot \bm{u}_h\right) \mathbb{I} \right):\nabla_h \ubold_h}{T_h},
 \end{align}
 and
 \begin{align}
 \frac{1}{C_v} \ipt{\left(\nu - \nu_h^{\mathrm{turb}}\right) \left(\nabla_h \bm{u}_h + \nabla_h \bm{u}_{h}^{T} - \frac{2}{3}\left(\nabla \cdot \bm{u}_h\right) \mathbb{I} \right):\nabla_h \ubold_h}{T_h} \leq 0,
 \end{align}
 when $\nu_h^{\mathrm{turb}} \geq \nu$.
\end{enumerate}
Under the assumptions outlined above, the viscous dissipation term cannot negatively impact the inf-sup stability of the associated method. Therefore, subject to these assumptions and the analysis in~\cite{miller2022versatile}, we can conclude that the finite element methods in Eqs.~\eqref{mass_cons_disc}--\eqref{energy_cons_disc} are inf-sup stable in the incompressible limit, and thus, are incompressibly stable. 


Lastly, we note that our versatile methods maintain stability of the discrete kinetic energy field (for incompressible flows with negligible viscous dissipation). In particular, if we construct a versatile method with BDM elements, we recover the kinetic-energy-stabilized H(div)-conforming method of~Guzman et al.~\cite{guzman2016h}.

\section{Numerical Experiments}

In this section, the results from several numerical simulations are presented for the finite element methods in Eqs.~\eqref{mass_cons_disc}--\eqref{energy_cons_disc}. These simulations were performed using Taylor-Hood elements with polynomials of degree $k$, $k$, and $k+1$ for the density, temperature and velocity respectively, unless otherwise stated. In the following simulations, the convective numerical fluxes were computed using upwind biased parameters $\zeta = \delta = 0.5$, and the viscous numerical fluxes were computed with $\eta = \varepsilon = 3(k+1)(k+2)$. The modeling constants $C_{\mathrm{mod}} = 0$ and $C_{\mathrm{smag}} = 0$ were used. In addition, a high-order BDF5 scheme was used for the time discretization. The meshes for each simulation were generated with quadrilateral elements split along the diagonal to make triangular elements. In the following sections, mesh dimensions will be reported as $N \times M$, where $N$ and $M$ refer to the number of quadrilaterals in the $x$ and $y$ directions. The actual total number of triangular elements is then $2(N \times M)$. Finally, all computations were carried out using the open-source finite element software package FEniCS~\cite{alnaes2015fenics}. FEniCS contains conforming and non-conforming finite element spaces, Newton’s method for solving nonlinear systems, and various Krylov iterative methods for solving linear systems. In this work, we have elected to use FEniCS’s Newton’s method in conjunction with the generalized minimal residual method (GMRES) linear solver and a Jacobi preconditioner. We note that our implementation is not optimized, and was constructed to demonstrate feasibility rather than exceptional speed.

The rest of this section contains four test cases for assessing the proposed versatile finite element methods. In section~\ref{MMS}, we evaluate their order of accuracy using the method of manufactured solutions. In section~\ref{APE}, we examine the compressible solution behavior at low Mach numbers, and prove that it converges to the correct incompressible solution.
In section~\ref{Cyl}, we perform a Mach number sweep on a two-dimensional cylinder in cross flow. Finally, in section~\ref{JAF}, we examine the drag forces on a Joukowski airfoil. 

\subsection{Method of Manufactured Solutions}\label{MMS}

For the first test case, we compared the results of our versatile method equipped with Taylor-Hood elements to an exact manufactured  solution. We defined the exact solution for $t \geq 0$ as follows
\begin{align*}
\rho &= \sin({x})\sin({y}) \exp(-2\nu t), \\[1.5ex]
T &= \frac{1}{2}\sin({x})\sin({y}) \exp\left(-2\frac{\kappa}{C_v}t\right),\\[1.5ex]
u &= \sin({x})\cos({y}) \exp(-2\nu t), \\[1.5ex]
v &= -\sin({y})\cos({x}) \exp(-2\nu t),
\end{align*}
where $\ubold = \left(u, v\right)$. Here, we used $\Omega = [0,1.25]^2$ as the spatial domain for our solution.  In addition, we note that the exact solution (above) was ``manufactured" by using the following forcing functions
\begin{align*}
S_{\rho} &= -\frac{2 \nu \sin({x})\sin({y})}{ \exp(2\nu t)}, \\[1.5ex]
S_T &= \Bigg(-4\nu\cos^{2}({x})\cos^{2}({y}) + \kappa\exp\left(2\left(\nu-\frac{\kappa}{C_v}\right)t \right)\sin({x})\sin({y}) \exp(2\nu t) \nonumber \\
&- (\kappa + C_v \nu)\sin({x})\sin({y})\Bigg)/ \exp(4 \nu t), \\[1.5ex]
 S_{u} &=\Bigg(\sin({x})\Big(\cos({x})\cos^2({y})\sin({x})\sin({y})\nonumber\\
 &+2 \mu \exp(2\nu t) \cos({y}) \exp(2\nu t) - 2\sin({x})\sin({y})\nonumber\\
 &+\cos({x})\sin^2({y}) ( \exp\left(\frac{-2 \kappa t}{C_v}+4\nu t\right)R + \sin({x})\sin({y}))\Big)\Bigg)/ \exp(6\nu t), \\[1.5ex]
S_{v} &=\Bigg(\sin({y})\Big(\cos({y})\cos^2({x})\sin({x})\sin({y})\nonumber \\
&-2 \mu \exp(2\nu t) \cos({x}) \exp(2\nu t)-2\sin({x})\sin({y})\nonumber \\
&+\cos({y})\sin^2({x}) ( \exp\left(\frac{-2 \kappa t}{C_v}+4\nu t\right)R + \sin({x})\sin({y}))\Big)\Bigg)/ \exp(6\nu t),
\end{align*}
where $\bm{S}_{u} = (S_u, S_v)$. We ran the simulation for $t \in [0,0.25]$, with parameters $C_v = 1$, $R =1$, $\nu = 3$, and $\kappa = 0.47$. The simulations were performed on a uniform grid with the exact solution specified as a Dirichlet boundary condition on all boundaries. Time stepping was performed with a step-size of $\Delta t = 5 \times 10^{-4}$. In addition, we utilized polynomials of degrees $k = 1$ and $2$. We therefore expected to see convergence rates of $k+2$ for the velocity and $k+1$ for the temperature and density, as the latter quantities are represented with a polynomial one-degree lower than the velocity. From table~\ref{MMS_tab}, we can see that the expected rates of convergence are achieved. 
\begin{table}[h!]
\begin{tabular}{|l|l|l||l|l||l|l||l|l|}
\hline
\multirow{2}{*}{k} & \multirow{2}{*}{h} & \multirow{2}{*}{dofs} & \multicolumn{2}{l||}{Velocity}        & \multicolumn{2}{l||}{Density}        & \multicolumn{2}{l|}{Temperature}     \\ \cline{4-9}
                   &                    &                       & $L^2$ error & order   & $L^2$ error & order   & $L^2$ error & order   \\ \hline
\multirow{4}{*}{1}  & 0.44194           & 212                  & 8.292e-5                  & -       & 0.001756                     & -       & 0.001954                    & -       \\ \cline{2-9}
                   & 0.22097           & 740                  & 1.047e-5                 & 2.984 & 4.470e-4                    & 1.974  & 4.974e-4                 & 1.974 \\ \cline{2-9}
                   & 0.11048            & 2756                 & 1.312e-6                  & 2.996 & 1.122e-4                   & 1.993 & 1.249e-4                  & 1.993  \\ \cline{2-9}
                   & 0.05524           & 10628                & 1.641e-7                  & 2.999 & 2.809e-5                   & 1.998 & 3.126e-5                 & 1.998 \\ \hline
\multirow{4}{*}{2} & 0.44194            & 500                   & 3.830e-6                   & -       & 9.1878e-5                    & -       & 1.022e-4                     & -       \\ \cline{2-9}
                   & 0.22097            & 1828                   & 2.392e-7                    & 4.001 & 1.148e-5                    & 2.999 & 1.278e-5                   & 2.999 \\ \cline{2-9}
                   & 0.11048            & 6980                  & 1.494e-8                   & 4.000 & 1.435-6                  & 2.999 & 1.597e-6                  & 2.999 \\ \cline{2-9}
                   & 0.05524             & 27268                 & 9.341e-10                  & 4.000 & 1.794e-7                   & 2.999  & 1.997e-7                 & 2.999 \\ \hline
\end{tabular}
\caption{Velocity, density, and temperature $L^2$ errors for various polynomial degrees $k$ and maximum element diameters $h$, for the versatile mixed method with Taylor-Hood elements. }
\label{MMS_tab}
\end{table}

\subsection{Asymptotic Preservation}\label{APE}

The second test case involved an isothermal vortex in a box with isothermal no-slip boundary conditions on each wall. The domain was $\Omega \in [0,1] \times [0,1]$, and we  tessellated it with a $64 \times 64$ triangular mesh. This single mesh was used in all of our numerical experiments. The case was run over the time interval $t \in [0,0.2]$. The fluid properties were set to $\mu = 0.01$, $C_v = 717.8$, $R = 287$, and $\gamma = 1.4$.

The goal of the test case was to show that the compressible simulation solutions converge to the incompressible solution in the presence of decreasing Mach number. In order to facilitate this comparison, we prescribed an initial, divergence-free, velocity field for both cases. The compressible case had the density and velocity specified, while the incompressible case had the velocity and kinematic pressure specified. The initial condition for the compressible case was specified as 
\begin{align*}
    \rho &= \rho_{ref} - \frac{\text{Ma}^2}{2}  \tanh(y - 0.5), \\[1.5ex]
    u &= \sin^2(\pi x) \sin(2 \pi y), \\[1.5ex]
    v &=-\sin^2(\pi  y)\sin(2\pi x).
\end{align*}
The initial condition for the incompressible case used the same velocity, but the kinematic pressure was specified as  
\begin{align*}
    p =(\rho_{ref})^{\gamma}.
\end{align*}
The incompressible case was simulated with constant density $\rho_{ref} = 1$, and the compressible case was simulated with Ma = 0.1, 0.05, 0.01, 0.005, 0.001, 0.0005, and 0.0001. The incompressible simulation was performed with Taylor-Hood elements of degree $k = 2$, and the compressible simulations were performed with two different methods: i) a versatile mixed method with Taylor-Hood elements of degree $k = 2$, and ii) a versatile mixed method with BDM elements of degree $k=2$. 

The time-steps for each compressible simulation are summarized in table~\ref{times_table}, and the time-step for the incompressible simulation was $5 \times 10^{-6}$. The $\Delta t$'s were chosen carefully in order to control the amount of temporal error that was generated, and to maintain the stability of each finite element scheme in conjunction with Newton's method.
\begin{table}[h!]
\centering 
\begin{tabular}{|l||l||l|}
\hline
Ma & TH, $\Delta t$     & BDM, $\Delta t$  \\ 
\hline
$0.1$ & 2.0e-6 & 2.0e-6  \\ \hline
$0.05$ & 2.0e-6 & 2.0e-6  \\ \hline
$0.01$  & 2.0e-6 & 2.0e-6 \\ \hline
$0.005$  & 2.0e-6  & 2.0e-6 \\ \hline
$0.001$  & 1.0e-6  & 5.0e-7 \\ \hline
$0.0005$  & 5.0e-7  & 1.0e-7 \\ \hline
$0.0001$  & 1.0e-7  & 5.0e-8 \\ \hline
\end{tabular}
\caption{Time-steps $\Delta t$ for two different finite element methods: a) the versatile mixed method with Taylor-Hood elements of degree $k=2$, and b) the versatile mixed method with BDM elements of degree $k=2$.}
\label{times_table}
\end{table}

At each Mach number, we calculated the differences between the incompressible and compressible approximations for the kinematic pressure and density. In particular, we calculated
\begin{align*}
    \left\| p_{\mathrm{comp}} - p_{\star} \right\|_{L^{2}(\Omega)}, \qquad \left\| \rho_{\mathrm{comp}} - \rho_{\star} \right\|_{L^{2}(\Omega)},
\end{align*}
where $\rho_{\mathrm{comp}}$ is the density extracted from the compressible simulations at different Mach numbers,
\begin{align*}
    p_{\mathrm{comp}} = (\rho_{\mathrm{comp}})^{\gamma},
\end{align*}
and $p_{\star}$ and $\rho_{\star}$ are computed as follows
\begin{align*}
    p_{\star} = 1 + \mathrm{Ma}^{2} p_{\mathrm{incomp}}, \qquad \rho_{\star} = \left(1 + \mathrm{Ma}^{2} p_{\mathrm{incomp}} \right)^{1/\gamma}.
\end{align*}
The quantity $p_{\mathrm{incomp}}$ is the kinematic pressure which was extracted from the incompressible simulation. Evidently, as Ma $\rightarrow 0$, we expect
\begin{align*}
    \left\| p_{\mathrm{comp}} - p_{\star} \right\|_{L^{2}(\Omega)} \longrightarrow 0, \qquad \left\| \rho_{\mathrm{comp}} - \rho_{\star} \right\|_{L^{2}(\Omega)} \longrightarrow 0.
\end{align*}

\begin{table}[h!]
\centering 
\begin{tabular}{|l||l|l||l|l|}
\hline
\multirow{2}{*}{Ma} & \multicolumn{2}{l||}{TH, $L^2$-differences}     & \multicolumn{2}{l|}{BDM, $L^2$-differences}      \\ \cline{2-5}
                   & Pressure  & Density  & Pressure  & Density   \\ \hline
$0.1$ & 4.269e-3 & 3.052e-3 & 4.075e-3 & 2.913e-3  \\ \hline
$0.05$ & 1.232e-3 & 8.799e-4 & 7.938e-4  & 5.671e-4  \\ \hline
$0.01$  & 3.888e-5 & 2.777e-5 & 4.057e-5 & 2.898e-5  \\ \hline
$0.005$  & 9.471e-6  & 6.765e-6 & 8.380e-6 & 5.986e-6  \\ \hline
$0.001$  & 4.525e-7  & 3.232e-7 & 4.017e-7 & 2.870e-7  \\ \hline
$0.0005$  & 7.900e-8  & 5.643e-8 & 8.350e-8 & 5.965e-8  \\ \hline
$0.0001$  & 5.155e-9  & 3.682e-9 & 5.360e-9 & 3.828e-9  \\ \hline
\end{tabular}
\caption{$L^2$ norms of differences between compressible and incompressible field variables (kinematic pressure and density), for different Mach numbers, and two different finite element methods: a) the versatile mixed method with Taylor-Hood elements of degree $k=2$, and b) the versatile mixed method with BDM elements of degree $k=2$.}
\label{Incomp_error_tab}
\end{table}

\begin{figure}[h!]
    \centering
    \includegraphics[width=0.7\linewidth]{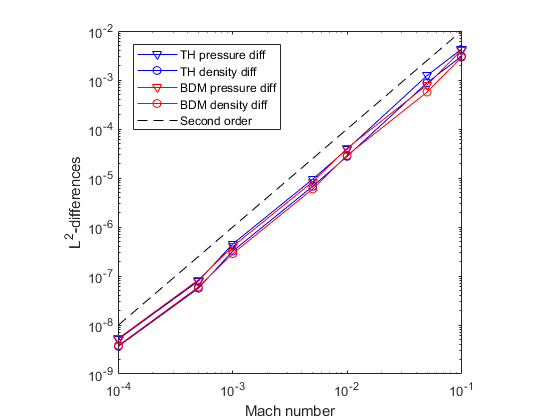}
    \caption{$L^2$ norms of differences between compressible and incompressible field variables (kinematic pressure and density), for two different finite element methods: a) the versatile mixed method with Taylor-Hood elements of degree $k=2$, and b) the versatile mixed method with BDM elements of degree $k=2$. The differences are shown to converge at a rate of 2nd order.}
    \label{fig:AP_converge}
\end{figure}

\begin{figure}[h!]
    \centering
    \includegraphics[width=0.55\linewidth]{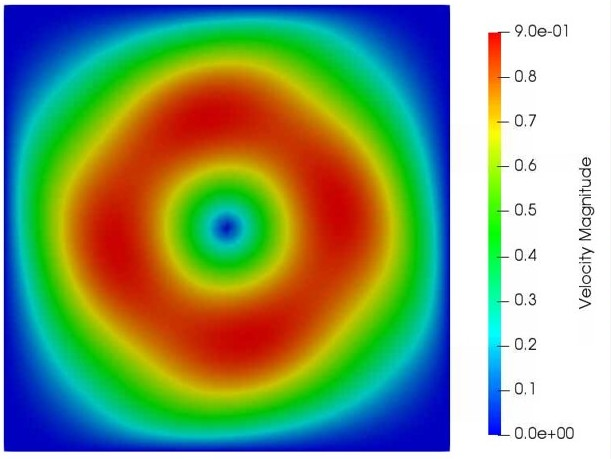}
    \caption{Velocity magnitude contours for the versatile mixed method with BDM elements of degree $k = 2$ at Ma = $0.005$ and $t = 0.18$.}
    \label{fig:AP_comp}
\end{figure}

We can see in table~\ref{Incomp_error_tab} that our mixed methods follow the expected trend at low Mach numbers. Moreover, we obtain 2nd-order convergence of the $L^2$-differences in kinematic pressure and density, as shown in Figure~\ref{fig:AP_converge}.

In addition, Figure~\ref{fig:AP_comp} shows that our BDM method maintains the initial decaying vortex. 

We believe that these results highlight a clear advantage of our approach: namely, since we are using function spaces traditionally associated with incompressible flows, we are successfully able to capture near-incompressible flow, even when solving the compressible equations.

\subsection{2-D Cylinder in Cross Flow}\label{Cyl}

For the third test case, we simulated a two-dimensional cylinder in cross flow over a range of Mach numbers, at a fixed Reynolds number Re = $100$.  In particular, the simulations were run with a uniform inlet flow, where the Mach number was adjusted over the following range of values: Ma = $0.1$, $0.2$, $0.3$, and $0.4$.  For these experiments, we utilized an initial density of $\rho$ = $1.0$ kg/m$^3$. The viscosity inside of the domain was allowed to change in accordance with Sutherland's law. The fixed Reynolds number at the inlet was maintained by adjusting the inlet viscosity in accordance with the inlet Mach number. In addition, a fixed  Prandtl number Pr = $0.72$ was used. The following fluid properties were assumed: $C_v$ = $717.8$ J/kg-K and $\gamma$ = $1.4$. The boundary condition for the left-most boundary was a subsonic inlet condition, and for the right-most boundary it was an extrapolation condition. The top and bottom boundaries utilized symmetry boundary conditions. The cylinder itself was equipped with adiabatic no-slip walls. 

The computational domain was $\Omega$ = $[-20,40]$ $ \times [-20,20]$. The cylinder had a diameter of $0.1247$ m and was centered at the point $(0,0)$. We used a mesh with~$70,000$ unstructured triangular elements. The polynomial order was $k =2$, and the simulations were run for $t \in [0,200]$ with a time-step of $\Delta t =$ $1.5 \times 10^{-5}$. Data for post-processing purposes was sampled over $t \in [150,200]$.

For this study, the primary quantity of interest was the time-averaged drag coefficient of the cylinder, $C_d$, defined as 
\begin{align}
   C_d = \frac{\overline{F}_d}{\frac{1}{2}\rho_\infty u_{\infty}^{2} d}. 
\end{align}
Here, $\overline{F}_d$ is the time-averaged drag force acting on the cylinder, $\rho_\infty$ is the free-stream density, $u_\infty$ is the free-stream flow speed, and $d$ is the diameter of the cylinder. 

\begin{figure}[h!]
    \centering
    \includegraphics[width=0.9\linewidth]{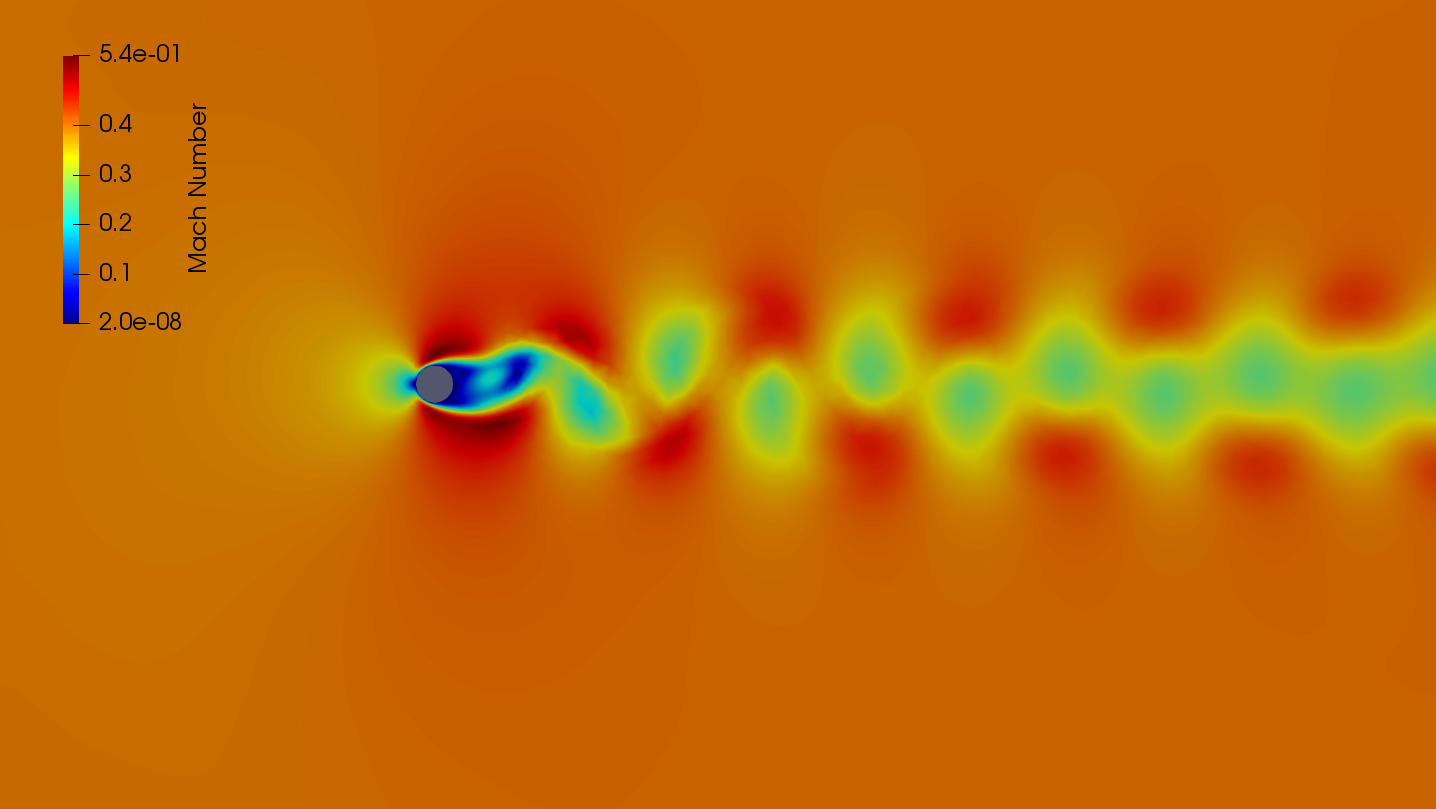}
    \caption{A snapshot of Mach number contours for the 2-D cylinder in cross flow at Ma = 0.4. Results were obtained using the versatile mixed method with Taylor-Hood elements of degree $k=2$.}
    \label{fig:Cyl_ma5}
\end{figure}

In what follows, we compare our results against an earlier study performed by~\cite{Canuto2015} on an identical geometry. At the given Reynolds number (100), the reference predicts that the flow will be unsteady with an oscillatory wake structure. This behavior is independent of the Mach number, as long as Ma $< 0.6$. As such, the flow fields for the various Mach numbers all exhibit similar behavior, which is correctly predicted for the Ma = 0.4 case by our method, (see figure ~\ref{fig:Cyl_ma5}). In particular, we can see from the figure the previously mentioned oscillatory wake behavior with the sinusoidal streamlines leaving the cylinder. This behavior is consistent across all Mach numbers, and implies that at least qualitatively we are in agreement with the reference solution. 
\begin{figure}[h!]
    \centering
    \includegraphics[width=0.9\linewidth]{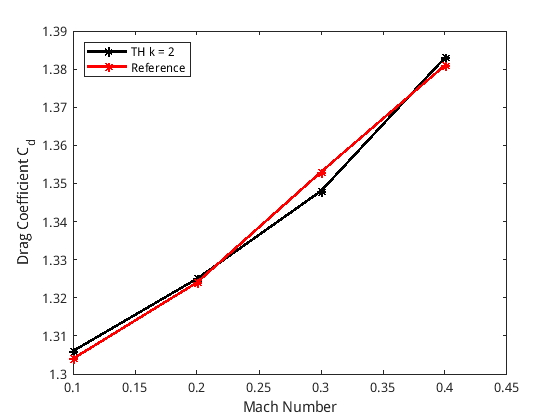}
    \caption{A plot of the time-averaged drag coefficient vs.~the Mach number. Results were obtained using the versatile mixed method with Taylor-Hood elements of degree $k=2$, and the reference~\cite{Canuto2015}.}
    \label{fig:Cyl_drag}
\end{figure}

In figure ~\ref{fig:Cyl_drag}, we compare our predictions for the drag coefficients with those of the reference solution~\cite{Canuto2015}. The predicted drag coefficient is in good agreement with the reference across the full range of Mach numbers considered.


\subsection{Joukowski Airfoil}\label{JAF}

The final test case involved flow over a Joukowski airfoil. This case was introduced at the 4th International Workshop on High-order CFD Methods,~(see \cite{iwohm4}). Here, an airfoil is simulated at an angle of attack of 0 degrees. The flow has Reynolds and Mach numbers of Re = 1000 and Ma = 0.5. The Prandtl number is again fixed at Pr = 0.72. Fluid properties for all simulations were $C_v = 717.8$ J/kg-K, $\mu = 1.716 \times 10^{-5}$ kg/m-s, and $\gamma = 1.4$. The viscosity was varied via Sutherland's Law. The initial density was prescribed as $\rho$ = 1 kg/m$^3$. The boundary conditions were specified as a subsonic inlet condition on the left-most boundary,  and an extrapolation condition on the right-most boundary. The airfoil itself was equipped with an adiabatic, no-slip condition.

The domain was $\Omega$ = $[-100,100]$ $\times [-100, 100]$ with the airfoil starting at $(0,0)$. The length of the airfoil from nose to trailing edge was 1 m. Unstructured meshes with 16,384, 65,536, and 262,144 triangular elements were used to tessellate the domain. These meshes were provided by organizers of the 5th International Workshop on High-order CFD Methods, and are numbered as meshes 2, 3, and 4, (see~\cite{iwohm5}). The simulations were run with polynomial order $k=2$ on these meshes until the drag converged to a steady-state value.
\begin{figure}[h!]
    \centering
    \includegraphics[width=0.9\linewidth]{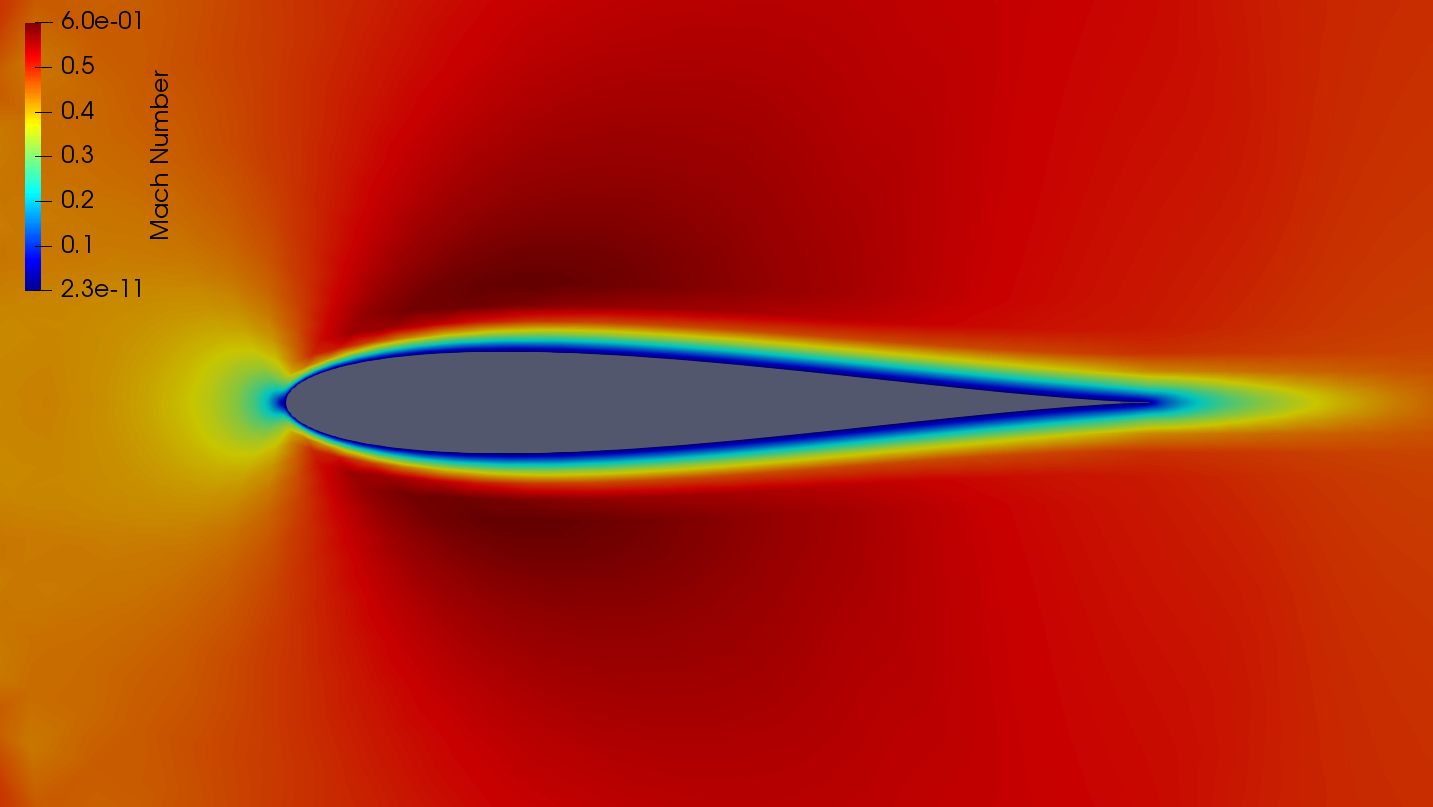}
    \caption{Mach number contours for the Joukowski airfoil at Re = $1000$ and Ma = $0.5$. Results were obtained using the versatile mixed method with Taylor-Hood elements of degree~$k=2$ on the finest mesh.}
    \label{fig:JAF}
\end{figure}
\begin{table}[h!]
\centering
\setlength{\extrarowheight}{2pt}
\begin{tabular}{|l||l|l|}
\hline
Mesh No. &DoF Count   & $C_d$   \\ \hline
2 & 16,704  & 0.1320                        \\ \hline
3 & 66,176  & 0.1273                     \\ \hline
4 & 263,424  & 0.1221                   \\ \hline
\end{tabular}
\caption{A table of drag coefficients on the given meshes. Results were obtained using the versatile mixed method with Taylor-Hood elements of degree $k=2$. The reference value for this case is $C_d = 0.1219$}
\label{JAF_drag}
\end{table}

Our results on mesh 4 are shown in figure~\ref{fig:JAF}. Here, the flow around the airfoil appears to be laminar, and the wake coming off the trailing edge has reached a steady state. 

The primary challenge for this case is to converge to the reference drag coefficient using the meshes provided. In table~\ref{JAF_drag}, we see that, as we increase the degrees of freedom in the simulation,  we move closer to the reference steady-state value of $C_d = 0.1219$. On the final grid, there is reasonable agreement between our predicted value and the reference.

\section{Conclusion}

In this work, we introduce a new class of ``versatile mixed finite element methods'' for solving the compressible Navier-Stokes equations, at low to moderate Mach numbers. These methods appear to be unique, as they simultaneously leverage multiple stabilization strategies for solving problems in the incompressible and compressible flow regimes. More precisely, our methods leverage numerical-flux-based stabilization, kinetic-energy-based stabilization, and inf-sup-based stabilization strategies. 

We note that our philosophy for designing these methods is somewhat unusual. In particular, many  numerical methods are only designed to be stable for linear advection and diffusion problems. In contrast, we have created finite element methods which are provably stable for the non-linear, non-isothermal, incompressible Navier-Stokes equations. Therefore, the starting point for our methods is significantly more complex than most, which facilitates their robustness and flexibility, in terms of successful application to increasingly complex problems. This claim has been demonstrated within the present paper, where we have successfully extended our methods to solve weakly-compressible flows.

The proposed methods are primarily designed to maintain their performance in the incompressible limit. In fact, we have ensured that the methods are `incompressibly stable', which means they are mathematically guaranteed to maintain stability when the density, dynamic viscosity coefficient, and heat conductivity coefficient assume constant values, and the viscous dissipation becomes negligible. We have shown through numerical experiments, that the resulting methods maintain good convergence properties and stability as the Mach number approaches zero. Based on this evidence, we argue that most (if not all) methods designed for weakly-compressible flows should satisfy the `incompressible stability' property. 

Lastly, we note that the new methods perform well, even far away from the incompressible limit, for flows in which the local Mach number exceeds 0.5. 
Of course, our numerical experiments in this regard are not exhaustive. Therefore, subsequent work will be needed to investigate the properties of these methods at higher Mach numbers (and Reynolds numbers), in an effort to establish their validity in these more challenging contexts, especially those involving shockwaves. 

\section*{Declaration of Competing Interests}

The authors declare that they have no known competing financial interests or personal relationships that could have appeared to influence the work reported in this paper.

\section*{Funding}

This research did not receive any specific grant from funding agencies in the public, commercial, or not-for-profit sectors.


\appendix

\section{Derivation of the Finite Element Methods} \label{method_deriv}

\subsection{Mass Equation Derivation}

%
One may substitute $\rho_h$ and $\ubold_h$ into Eq.~\eqref{mass_cons}, multiply by a test function $q_h$, and integrate over the entire domain in order to yield
\begin{align}
\ipt{\partial_t \, \rho_h}{q_h} + \ipt{\nabla_h \cdot \left( \rho_h \ubold_h\right)}{q_h} = \ipt{S_{\rho}}{q_h}. \label{mass_ibp_one}
\end{align}
This is identical to Eq.~\eqref{mass_cons_disc}.
%


\subsection{Linear Momentum Equation Derivation}

One may substitute $\rho_h$, $\ubold_h$, and $T_h$ into Eq.~\eqref{moment_cons}, compute the dot product with a test function $\wbold_h$, and integrate over the entire domain in order to yield
\begin{align}
\ipt{\partial_t \left( \rho_h \ubold_h\right)}{\wbold_h} + \ipt{\nabla_h \cdot \left( \left(\rho_h \ubold_h \right) \otimes \ubold_h + P_h \mathbb{I} \right)}{\wbold_h} - \ipt{ \nabla_h \cdot \left( \rho_h \bm{\tau}_h \right)}{\wbold_h} = \ipt{\bm{S}_{u}}{\wbold_h}. \label{moment_ibp_neg}
\end{align}
Upon integrating the second and third terms on the LHS by parts and inserting numerical fluxes $\widehat{\bm{\sigma}}_{\text{inv}}$ and $\widehat{\bm{\sigma}}_{\text{vis}}$, one obtains
\begin{align}
&\ipt{\nabla_h \cdot \left( \left(\rho_h \ubold_h \right) \otimes \ubold_h + P_h \mathbb{I} \right)}{\wbold_h} \label{moment_ibp_zero} \\[1.5ex]
\nonumber &= - \ipt{ \left(\rho_h \ubold_h \right) \otimes \ubold_h + P_h \mathbb{I}}{\nabla_h \wbold_h} + \ipbt{\left( \left(\rho_h \ubold_h \right) \otimes \ubold_h + P_h \mathbb{I} \right)\nbold}{\wbold_h}  \\[1.5ex]
\nonumber & \equiv - \ipt{ \left(\rho_h \ubold_h\right) \otimes \ubold_h + P_h \mathbb{I}}{\nabla_h \wbold_h} + \ipbt{ \widehat{\bm{\sigma}}_{\text{inv}} \, \nbold}{\wbold_h} \\[1.5ex]
\nonumber & = - \ipt{ \left(\rho_h \ubold_h\right) \otimes \ubold_h}{\nabla_h \wbold_h} - \ipt{P_h}{\nabla \cdot \wbold_h} + \ipbt{ \widehat{\bm{\sigma}}_{\text{inv}} \, \nbold}{\wbold_h}.
\end{align}
\begin{align}
-\ipt{ \nabla_h \cdot \left( \rho_h \bm{\tau}_h\right)}{\wbold_h} &= \ipt{\rho_h \bm{\tau}_h}{\nabla_h \wbold_h} - \ipbt{\rho_h \bm{\tau}_h \, \nbold}{\wbold_h} \label{moment_ibp_one} \\[1.5ex]
\nonumber & \equiv \ipt{\rho_h \bm{\tau}_h}{\nabla_h \wbold_h} - \ipbt{ \widehat{\bm{\sigma}}_{\text{vis}} \, \nbold}{\wbold_h}. 
\end{align}
Consider substituting the definition of $\bm{\tau}_h$ (Eq.~\eqref{stress_tensor}) into the first term on the RHS of Eq.~\eqref{moment_ibp_one}
\begin{align}
\ipt{\rho_h \bm{\tau}_h}{\nabla_h \wbold_h} = \ipt{\mu_h \left( \nabla_h \ubold_h + \nabla_h \ubold_h^T - \frac{2}{3} \left(\nabla \cdot \ubold_h \right) \mathbb{I} \right)}{\nabla_h \wbold_h}. \label{moment_ibp_two}
\end{align}
One may expand each term in Eq.~\eqref{moment_ibp_two} by integrating by parts, inserting a numerical flux $\widehat{\bm{\varphi}}_{\text{vis}}$, and integrating by parts again as follows
\begin{align}
\ipt{\mu_h \nabla_h \ubold_h}{\nabla_h \wbold_h} & = \ipt{\nabla_h \ubold_h}{\mu_h \nabla_h \wbold_h} \label{moment_ibp_three} \\[1.5ex]
\nonumber &= -\ipt{\ubold_h}{\nabla_h \cdot \left(\mu_h \nabla_h \wbold_h \right)} + \ipbt{\mu_h \ubold_h}{ \left(\nabla_h \wbold_h \right) \nbold} \\[1.5ex]
\nonumber & \equiv -\ipt{\ubold_h}{\nabla_h \cdot \left(\mu_h \nabla_h \wbold_h \right)} + \ipbt{\widehat{\bm{\varphi}}_{\text{vis}}}{ \left(\nabla_h \wbold_h \right) \nbold} \\[1.5ex]
\nonumber & = \ipt{\mu_h \nabla_h \ubold_h}{\nabla_h \wbold_h} + \ipbt{\widehat{\bm{\varphi}}_{\text{vis}} - \mu_h \ubold_h }{ \left(\nabla_h \wbold_h \right) \nbold}.
\end{align}
\begin{align}
\ipt{\mu_h \nabla_h \ubold_h^T}{\nabla_h \wbold_h} & = \ipt{\nabla_h \ubold_h}{\mu_h \nabla_h \wbold_h^T} \label{moment_ibp_four} \\[1.5ex]
\nonumber &= -\ipt{\ubold_h}{\nabla_h \cdot \left(\mu_h \nabla_h \wbold_h^T \right)} + \ipbt{\mu_h \ubold_h}{ \left(\nabla_h \wbold_h^T \right) \nbold} \\[1.5ex]
\nonumber & \equiv -\ipt{\ubold_h}{\nabla_h \cdot \left(\mu_h \nabla_h \wbold_h^T \right)} + \ipbt{\widehat{\bm{\varphi}}_{\text{vis}}}{ \left(\nabla_h \wbold_h^T \right) \nbold} \\[1.5ex]
\nonumber & = \ipt{\mu_h \nabla_h \ubold_h^T}{\nabla_h \wbold_h} + \ipbt{\widehat{\bm{\varphi}}_{\text{vis}} - \mu_h \ubold_h }{ \left(\nabla_h \wbold_h^T \right) \nbold}.
\end{align}
\begin{align}
\ipt{- \frac{2}{3} \mu_h \left(\nabla \cdot \ubold_h \right) \mathbb{I}}{\nabla_h \wbold_h} &= - \frac{2}{3} \ipt{\left(\nabla \cdot \ubold_h \right) }{\mu_h \nabla \cdot \wbold_h} \label{moment_ibp_five} \\[1.5ex]
\nonumber & = -\frac{2}{3} \left( -\ipt{\ubold_h}{ \nabla_h \left(\mu_h \nabla_h \cdot \wbold_h \right)} + \ipbt{\mu_h \ubold_h}{ \left(\nabla \cdot \wbold_h \right) \nbold} \right) \\[1.5ex]
\nonumber & \equiv -\frac{2}{3} \left( -\ipt{\ubold_h}{ \nabla_h \left(\mu_h \nabla_h \cdot \wbold_h \right)} + \ipbt{\widehat{\bm{\varphi}}_{\text{vis}}}{ \left(\nabla \cdot \wbold_h \right) \nbold} \right) \\[1.5ex]
\nonumber & = \ipt{- \frac{2}{3} \mu_h \left(\nabla \cdot \ubold_h \right) \mathbb{I}}{\nabla_h \wbold_h} -\frac{2}{3} \ipbt{\widehat{\bm{\varphi}}_{\text{vis}} - \mu_h \ubold_h}{ \left(\nabla \cdot \wbold_h \right) \nbold}.
\end{align}
Upon combining Eqs.~\eqref{moment_ibp_three} -- \eqref{moment_ibp_five} along with the definition of $\bm{\tau}_h$ (Eq.~\eqref{stress_tensor}), one obtains
\begin{align}
\ipt{\rho_h \bm{\tau}_h}{\nabla_h \wbold_h} \equiv \ipt{\rho_h \bm{\tau}_h}{\nabla_h \wbold_h} +  \ipbt{\widehat{\bm{\varphi}}_{\text{vis}} - \mu_h \ubold_h }{ \left(\nabla_h \wbold_h + \nabla_h \wbold_h^T -\frac{2}{3} \left( \nabla \cdot \wbold_h  \right) \mathbb{I} \right) \nbold}. \label{moment_ibp_six}
\end{align}
Finally, one may substitute Eqs.~\eqref{moment_ibp_zero}, \eqref{moment_ibp_one}, and \eqref{moment_ibp_six}, into Eq.~\eqref{moment_ibp_neg}, in order to obtain Eq.~\eqref{moment_cons_disc}.

\subsection{Internal Energy Equation Derivation}

One may substitute $\rho_h$, $\ubold_h$, and $T_h$ into Eq.~\eqref{energy_cons}, multiply by a test function $r_h$, and integrate over the entire domain in order to yield
\begin{align}
\nonumber &\ipt{\partial_t \left( \rho_h T_h \right)}{r_h}  + \ipt{\nabla_h \cdot \left( \rho_h T_h \ubold_h \right)}{r_h} - \ipt{\nabla_h \cdot \left( \frac{\kappa_h}{C_v} \nabla_h T_h  \right)}{r_h} \\[1.5ex]
&= - \left(\gamma-1\right) \ipt{\rho_h T_h \left(\nabla \cdot \ubold_h \right)}{r_h} + \frac{1}{C_v} \ipt{\rho_h \bm{\tau_h}:\nabla_h \ubold_h}{r_h} + \ipt{S_{T}}{r_h}. \label{energy_ibp_neg}
\end{align}
Upon integrating the second and third terms on the LHS by parts and inserting numerical fluxes $\widehat{\bm{\phi}}_{\text{inv}}$ and $\widehat{\bm{\phi}}_{\text{vis}}$, one obtains
\begin{align}
\ipt{\nabla_h \cdot \left( \rho_h T_h \ubold_h \right)}{r_h} &= -\ipt{\rho_h T_h \ubold_h}{\nabla_h r_h} + \ipbt{ \left( \rho_h T_h \ubold_h \right) \cdot \nbold}{r_h} \label{energy_ibp_zero} \\[1.5ex]
\nonumber & \equiv  -\ipt{\rho_h T_h \ubold_h}{\nabla_h r_h} + \ipbt{ \widehat{\bm{\phi}}_{\text{inv}} \cdot \nbold}{r_h}.
\end{align}
\begin{align}
 -\ipt{\nabla_h \cdot \left( \frac{\kappa_h}{C_v} \nabla_h T_h \right)}{r_h} &= \ipt{\frac{\kappa_h}{C_v} \nabla_h T_h}{\nabla_h r_h} - \ipbt{ \left(\frac{\kappa_h}{C_v} \nabla_h T_h \right) \cdot \nbold}{r_h} \label{energy_ibp_one} \\[1.5ex]
\nonumber & \equiv \ipt{\frac{\kappa_h}{C_v} \nabla_h T_h}{\nabla_h r_h} - \ipbt{ \widehat{\bm{\phi}}_{\text{vis}} \cdot \nbold}{r_h}.
\end{align}
One may further expand Eq.~\eqref{energy_ibp_one}. Consider  integrating by parts, inserting the numerical flux $\widehat{\lambda}_{\text{vis}}$, and integrating by parts again as follows
\begin{align}
\ipt{\nabla_h T_h}{\frac{\kappa_h}{C_v} \nabla_h r_h} & = - \ipt{T_h}{\nabla_h \cdot \left( \frac{\kappa_h}{C_v} \nabla_h r_h \right)} + \ipbt{\frac{\kappa_h}{C_v} T_h}{\nabla_h r_h \cdot \nbold} \label{energy_ibp_seven} \\[1.5ex]
\nonumber & \equiv - \ipt{T_h}{\nabla_h \cdot \left( \frac{\kappa_h}{C_v} \nabla_h r_h \right)} + \ipbt{\widehat{\lambda}_{\text{vis}}}{\nabla_h r_h \cdot \nbold} \\[1.5ex]
\nonumber & = \ipt{\frac{\kappa_h}{C_v} \nabla_h T_h}{\nabla_h r_h} + \ipbt{\widehat{\lambda}_{\text{vis}} - \frac{\kappa_h}{C_v} T_h}{\nabla_h r_h \cdot \nbold}.
\end{align}
Finally, one may combine Eqs.~\eqref{energy_ibp_neg}--\eqref{energy_ibp_seven} in order to obtain Eq.~\eqref{energy_cons_disc}.


{\scriptsize\bibliography{technical-refs}}

\end{document}